\documentclass[12 pt]{amsart}
\usepackage{amsmath,times,epsfig,amssymb,amsbsy,amscd,amsfonts,amstext,graphicx,color,bm,amsthm}
\usepackage[all]{xy}
\usepackage{graphicx}

\usepackage[urlcolor=blue]{hyperref}
\usepackage{tikz,subfigure}

\usepackage{listings} 

\newlength{\colheight}
\newlength{\colwidth}
\definecolor{main}{rgb}{0.3,0.7,0.9}
\definecolor{red-}{rgb}{1.0, 0.2, 0.0}
\definecolor{blue-}{rgb}{0.0, 0.2, 1.0}
\definecolor{green-}{rgb}{0.0, 0.6, 0.0}
\definecolor{greenblue}{rgb}{0.0, 0.3, 1}
\definecolor{gold}{rgb}{0.8, 0.7, 0.0}

\lstdefinelanguage{cocoa}
{
  basicstyle=\small\ttfamily,
  commentstyle=\color{red!90!black},       
  stringstyle=\color{green!70!black},      
  morecomment=[l]{//},
  morecomment=[l]{--},
  morecomment=[s]{/*}{*/},
  morestring=[b]{"}, 
  classoffset=1,
  morekeywords={
    define,enddefine,if,endif,for,endfor, 
    use,in,then,else,elif,return,
    and,or,
    break,continue,ciao,do,exit,
    ImportByValue,ImportByRef,importbyvalue,importbyref,
    in,isin,IsIn,
    on,opt,PrintLn,println,print,
    quit,ref,return,step,
    toplevel,TopLevel,
    then,to,
    use
  },
  keywordstyle=\color{blue!70!green!80!white},
  classoffset=2,
  morekeywords={
    Ideal,Mat,       Not,Record,Error,
    ideal,mat,matrix,not,record,error,submodule
  },
  keywordstyle=\color{purple!50!white},
  classoffset=3,
  morekeywords={
    TRUE,FALSE,True,False,true,false,
    Lex,Xel,DegLex,DegRevLex,
    PosTo,ToPos,Null
  },
  keywordstyle=\color{brown},
}

\newcommand{\C}{\mathbb{C}}

\newcommand{\R}{\mathbb{R}}

\renewcommand{\to}{\longrightarrow}

\newtheorem{Theorem}{Theorem}[section]
\newtheorem{Definition}[Theorem]{Definition}

\DeclareMathOperator{\Der}{Der}

\DeclareMathOperator{\pdeg}{pdeg}

\DeclareMathOperator{\codim}{codim}

\def \X(#1){\{x_1,\dots, x_{#1}\}}

\DeclareMathOperator{\rk}{rk}

\newcommand{\A}{\mathcal{A}}

\DeclareMathOperator{\OT}{OT}
\DeclareMathOperator{\AOT}{AOT}
\DeclareMathOperator{\ST}{ST}

\newcommand \ideal[1] {\langle #1 \rangle}

\begin{document}

\title{Hyperplane arrangements in CoCoA}

\begin{abstract}
We introduce the package \textbf{arrangements} for the software CoCoA. This package provides a data structure and the necessary methods for working with hyperplane arrangements. 
In particular, the package implements methods to enumerate many commonly studied classes of arrangements, perform operations on them, and calculate various invariants associated to them.
\end{abstract}

\author{Elisa Palezzato}
\address{Elisa Palezzato, Department of Mathematics, Hokkaido University, Kita 10, Nishi 8, Kita-Ku, Sapporo 060-0810, Japan.}
\email{palezzato@math.sci.hokudai.ac.jp}
\author{Michele Torielli}
\address{Michele Torielli, Department of Mathematics, Hokkaido University, Kita 10, Nishi 8, Kita-Ku, Sapporo 060-0810, Japan.}
\email{torielli@math.sci.hokudai.ac.jp}


\date{\today}
\maketitle


\section{Introduction}
An arrangement of hyperplanes is a finite collection of codimension one affine subspaces in a finite dimensional vector space. Associated to these spaces,
there is a plethora of algebraic, combinatorial and topological invariants. Arrangements are easily defined but they lead to
deep and beautiful results that put in connection various area of mathematics. We refer the reader to the seminal text \cite{orlterao} for a comprehensive account of this subject. 

One of the main goals in the study of hyperplane arrangements is to decide whether a given invariant is combinatorically determined, and, if so, to express it explicitly in terms of the intersection lattice of the arrangement.

We describe the new package \textbf{arrangements} that computes several combinatorial invariants
(like the lattice of intersections and its flats, the Poincar\'e, the characteristic and the Tutte polynomials)
 and algebraic ones (like the Orlik-Terao and the Solomon-Terao ideals)
of hyperplane arrangement for the software CoCoA (\cite{CoCoALib}, \cite{AbbottBigatti2016} and \cite{COCOA}). 
Moreover, several functions for the class of free hyperplane arrangements are implemented. In addition, this package allows 
also to do computations with multiarrangements.
Finally, several known families of arrangements (like classic reflection arrangements, Shi arrangements, Catalan arrangements, Shi-Catalan arrangements, graphical arrangements
and signed graphical ones) can be easily constructed.

We introduce this package via several examples.  
Specifically, in Section 2 we first recall the definitions of various combinatorial invariants of a given arrangement and then describe how to compute them. 
In Section 3, we describe how to work with free hyperplane arrangements, and in Section 4
how to define the Orlik-Terao ideal and the Solomon-Terao one. 
Finally, in Section 5 we describe the class of multiarrangements with particular emphasis to the free ones.

This package will be part of the official release CoCoA-5.2.4.

\section{Combinatoric of arrangements}

Let $V$ be a vector space of dimension $l$ over a field $K$. Fix a system of coordinate $(x_1,\dots, x_l)$ of $V^\ast$. 
We denote by $S = S(V^\ast) = K[x_1,\dots, x_l]$ the symmetric algebra. A finite set of affine hyperplanes $\A =\{H_1, \dots, H_n\}$ in $V$ is called a \textbf{hyperplane arrangement}. 

For each hyperplane $H_i$ we fix a defining equation $\alpha_i\in S$ such that $H_i = \alpha_i^{-1}(0)$, 
and let $Q(\A)=\prod_{i=1}^n\alpha_i$. An arrangement $\A$ is called \textbf{central} if each $H_i$ contains the origin of $V$. 
In this case, the defining equation $\alpha_i\in S$ is linear homogeneous, and hence $Q(\A)$ is a homogeneous polynomial of degree $n$. 

The operation of \textbf{coning} allows to transform any arrangement $\A$ of $V$ with $n$ hyperplanes
 into a central arrangement $c\A$ in a vector space of dimension $l+1$ with $n+1$ hyperplanes, see \cite{orlterao}.

Notice that in CoCoA to compute the cone of an arrangement $\A$, the homogenizing variable needs to be already present in the ring in which the equation of $\A$ is defined.
For example, we can construct the cone of the Shi arrangement of type $A$ in CoCoA as follows:
\begin{lstlisting}[alsolanguage=cocoa]
/**/ use S::=QQ[x,y,z,w];
/**/ A := ArrShiA(S, 3); A;
[x-y,  x-z,  y-z,  x-y-1,  x-z-1,  y-z-1]
/**/ ArrCone(A, w);
[x-y,  x-z,  y-z,  x-y-w,  x-z-w,  y-z-w,  w]
\end{lstlisting}

Let $L(\A)=\{\bigcap_{H\in\mathcal{B}}H \mid \mathcal{B}\subseteq\A\}$ be the \textbf{lattice of intersection} of $\A$. 
Define a partial order on $L(\A)$ by $X\le Y$ if and only if $Y\subseteq X$, for all $X,Y\in L(\A)$. 
Note that this is the reverse inclusion. The elements of $L(\A)$ are called \textbf{flats} of $\A$. Define a rank function on $L(\A)$ by $\rk(X)=\codim(X)$. 
$L(\A)$ plays a fundamental role in the study of hyperplane arrangements, in fact it determines the combinatoric of the arrangement.

We can compute the flats in the lattice of intersection of the reflection arrangement of type $D$ in CoCoA in the following way:

\begin{lstlisting}[alsolanguage=cocoa]
/**/ use S::=QQ[x,y,z];
/**/ A := ArrTypeD(S,3); A;
[x-y,  x+y,  x-z,  x+z,  y-z,  y+z]


/**/ ArrFlats(A);
[[ideal(0)], 
 [ideal(x-y), ideal(x+y), ideal(x-z), ideal(x+z), 
  ideal(y-z), ideal(y+z)], 
 [ideal(x, y), ideal(x-z, y-z), ideal(x+z, y+z), 
  ideal(x-z, y+z), ideal(x+z, y-z), ideal(x, z), 
  ideal(y, z)], 
 [ideal(x, y, z)]]
\end{lstlisting}

Let $\mu\colon L(\A)\to\mathbb{Z}$ be the \textbf{M\"obius function} of $L(\A)$ defined by
$$\mu(X)=
\begin{cases}
      1 & \text{for } X=V,\\
      -\sum_{Y<X}\mu(Y) & \text{if } X>V.
\end{cases}$$

The \textbf{Poincar\'e polynomial} of $\A$ is defined by $$\pi(\A,t) = \sum_{X\in L(\A)}\mu(X)(-t)^{\rk(X)},$$
and it satisfies the formula
$$\pi(c\A,t)=(t+1)\pi(\A,t).$$

We now verify the previous result for the Shi arrangement of type $A$ in CoCoA.
\begin{lstlisting}[alsolanguage=cocoa]
/**/ use S::=QQ[x,y,z,w];
/**/ A := ArrShiA(S, 3);
/**/ pi_A := ArrPoincarePoly(A); pi_A;
9*t^2 +6*t +1
/**/ cA := ArrCone(A, w);
/**/ pi_cA := ArrPoincarePoly(cA); pi_cA;
9*t^3+15*t^2+7*t+1
/**/ pi_A := ArrPoincarePoly(A);
/**/ t := indets(RingOf(pi_A),1);
/**/ pi_cA = (1+t)*pi_A;
true
\end{lstlisting}

For any flat $X\in L(\A)$ define the subarrangement $\A_X$ of $\A$ by 
$$\A_X=\{H\in\A~|~X\subseteq H\}.$$
Similarly, define the \textbf{restriction} of $\A$ to $X$ as the arrangement $\A^X$ in $X$
$$\A^X=\{X\cap H~|~H\in\A\setminus\A_X \text{ and } X\cap H\ne\emptyset\}.$$

The \textbf{characteristic polynomial} of $\A$ is $$\chi(\A,t) =t^l\pi(\A,-t^{-1})= \sum_{X\in L(\A)}\mu(X)t^{\dim(X)}.$$ 
The characteristic polynomial is characterized by the following recursive relation
$$\chi(\A, t) = \chi(\A_H, t)-\chi(\A^H, t),$$
for any $H\in \A$.

We verify the previous result for $\A^{[-1,2]}$ the Shi-Catalan arrangement of type $A$ in CoCoA.
\begin{lstlisting}[alsolanguage=cocoa]
/**/ use S ::= QQ[x,y,z];
/**/ A := ArrShiCatalanA(S, 3, [-1, 2]); A;
[x-y, x-z, y-z, x-y-1, x-z-1, y-z-1, x-y+1, x-y+2, 
 x-z+1, x-z+2, y-z+1, y-z+2]
/**/ A_1 := ArrDeletion(A,4); A_1;
[x-y, x-z, y-z, x-z-1, y-z-1, x-y+1, x-y+2, x-z+1, 
 x-z+2, y-z+1, y-z+2]
/**/ A_2 := ArrRestriction(A,4); A_2;
[y[1]-y[2]+1, y[1]-y[2], y[1]-y[2]-1, y[1]-y[2]+2, 
 y[1]-y[2]+3]
/**/ ArrCharPoly(A) = ArrCharPoly(A_1) - ArrCharPoly(A_2);
true
\end{lstlisting}

For $i=0,\dots,l$ we define the \textbf{$i$-th Betti number} $b_i(\A)$  by the formula
$$\chi(\A,t)=\sum_{i=0}^l(-1)^ib_i(\A)t^{l-i}.$$
The importance of the characteristic polynomial in combinatorics is justified by the following result 
from \cite{crapo1970foundations}, \cite{orlik1980combinatorics} and \cite{zaslavsky1975facing}.

\begin{Theorem} 
We have that
\begin{enumerate}
\item If $\A$ is an arrangement in $\mathbb{F}_q^l$ (vector space over a finite field $\mathbb{F}_q$), then $|\mathbb{F}^l_q\setminus \bigcup_{H\in\A}H|=\chi(\A, q)$.
\item If $\A$ is an arrangement in $\C^l$, then the topological $i$-th Betti number of the complement is $b_i(\C^l  \setminus \bigcup_{H\in\A}H)=b_i(\A)$.
\item If $\A$ is an arrangement in $\R^l$, then $|\chi(\A,-1)|$ is the number of chambers and $|\chi(\A, 1)|$ is the number of bounded chambers.
\end{enumerate}
\end{Theorem}

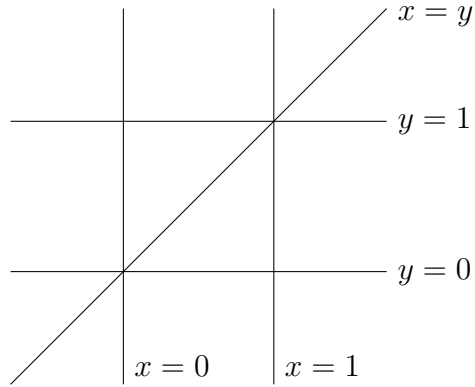
\begin{figure}[htbp]
\begin{tikzpicture}[scale=0.5] 
\draw(-5,-5) -- (5,5);
\draw(-2,-5) -- (-2,5);
\draw(2,-5) -- (2,5);
\draw(-5,-2) -- (5,-2);
\draw(-5,2) -- (5,2);

\coordinate [label=right:{$x=0$}] (V0) at (-2,-4.5);
\coordinate [label=right:{$x=1$}] (V0) at (2,-4.5);
\coordinate [label=right:{$y=0$}] (V0) at (5,-2);
\coordinate [label=right:{$y=1$}] (V0) at (5,2);
\coordinate [label=right:{$x=y$}] (V0) at (5,4.8);
\end{tikzpicture} 
\caption{A line arrangement in $\mathbb{R}^2$}\label{Fig:linearr}
\end{figure}

Using the previous statements, we can compute the Betti numbers, the number of chambers and the number of bounded chambers of the arrangement in Figure \ref{Fig:linearr} in CoCoA.
\begin{lstlisting}[alsolanguage=cocoa]
/**/ use S ::=QQ[x,y];
/**/ A := [x,x-1,y,y-1,x-y];
/**/ ArrBettiNumbers(A);
[1,  5,  6]

/**/ NumChambers(A); 
12
/**/ NumBChambers(A);
2
\end{lstlisting}

Associated to each hyperplane arrangement, it can be naturally defined a third polynomial.
The \textbf{Tutte polynomial} of $\A$ is
$$T_\A(x,y)=\sum_{\substack{\mathcal{B}\subseteq\A \\ \mathcal{B} \text{ central }}}(x-1)^{\rk(\A)-\rk(\mathcal{B})}(y-1)^{|\mathcal{B}|-\rk(\mathcal{B})}.$$
As shown in \cite{ardila2007computing}, it turns out that the Tutte and the characteristic polynomials are related by
$$\chi(\A,t)=(-1)^{\rk(\A)}t^{|\A|-\rk(\A)}T_{\A}(1-t,0).$$

We verify the previous result for the Boolean arrangement in CoCoA. Notice that here,
since the Tutte and the characteristic polynomials live in different rings, we need to construct a ring 
homomorphism, with the command \texttt{PolyRingHom}, to check the required equality.
\begin{lstlisting}[alsolanguage=cocoa]
/**/ use S ::= QQ[x,y,z];
/**/ A := ArrBoolean(S, 3);
/**/ Tutte_A := ArrTuttePoly(A); Tutte_A;
t[1]^3
/**/ char_A := ArrCharPoly(A);
/**/ R:=RingOf(Tutte_A);
/**/ P:=RingOf(char_A);
/**/ t:=indets(P,1);
/**/ psi := CanonicalHom(BaseRing(R),P);
/**/ phi := PolyRingHom(R, P, psi, [1-t,0]);
/**/ char_A=(-1)^3*t^(len(A)-3)*phi(Tutte_A);
true
\end{lstlisting}

\section{Free hyperplane arrangements}
In the theory of hyperplane arrangements, the freeness of an arrangement is a very important algebraic property. 
In fact, freeness implies several interesting geometric and combinatorial properties of the arrangement itself.
See for example \cite{terao1980arrangementsI}, \cite{yoshinaga2012freeness}, \cite{abe2016divisionally-im}, \cite{Gin-freearr} and \cite{palezzato2018free}.

We denote by $\Der_{V} =\{\sum_{i=1}^l f_i\partial_{x_i}~|~f_i\in S\}$ the $S$-module of \textbf{polynomial vector fields} on $V$ (or $S$-derivations). 
Let $\delta =  \sum_{i=1}^l f_i\partial_{x_i}\in \Der_{V}$. iI $f_1, \dots, f_l$ are homogeneous polynomials of degree~$d$ in $S$, then $\delta$ is  said to be \textbf{homogeneous of polynomial degree} $d$. 
In this case, we write $\pdeg(\delta) = d$.

For any central arrangement $\A$ we define the \textbf{module of vector fields logarithmic tangent} to $\A$ (logarithmic vector fields) by
$$D(\A) = \{\delta\in \Der_{V}~|~ \delta(\alpha_i) \in \ideal{\alpha_i} S, \forall i\}.$$

The module $D(\A)$ is obviously a graded $S$-module and we have that $D(\A)= \{\delta\in \Der_{V}~|~ \delta(Q(\A)) \in \ideal{ Q(\A)} S\}$. 

\begin{Definition} 
A central arrangement $\A$ is said to be \textbf{free with exponents $(e_1,\dots,e_l)$} 
if and only if $D(\A)$ is a free $S$-module and there exists a basis $\delta_1,\dots,\delta_l \in D(\A)$ 
such that $\pdeg(\delta_i) = e_i$, or equivalently $D(\A)\cong\bigoplus_{i=1}^lS(-e_i)$.
\end{Definition}



Let $\delta_1,\dots,\delta_l \in D(\A)$. 
Then $\det(\delta_i(x_j))$ is divisible by $Q(\A)$.
One of the most famous characterization of freeness is due to Saito \cite{saito} and it uses  the determinant of the coefficient 
matrix of $\delta_1,\dots,\delta_l$ to check if the arrangement $\A$ is free or not.

\begin{Theorem}[Saito's criterion]\label{theo:saitocrit}
Let $\delta_1, \dots, \delta_l \in D(\A)$. Then the following facts are equivalent
\begin{enumerate}
\item $D(\A)$ is free with basis $\delta_1, \dots, \delta_l$, i. e. $D(\A) = S\cdot\delta_1\oplus \cdots \oplus S \cdot\delta_l$.
\item $\det(\delta_i(x_j))=c Q(\A)$, where $c\in K\setminus\{0\}$.
\item $\delta_1, \dots, \delta_l$ are linearly independent over $S$ and $\sum_{i=1}^l\pdeg(\delta_i)=n$.
\end{enumerate}
\end{Theorem}
%

Given a simple graph $G$, we can define the graphical arrangement $\A(G)$, see \cite{orlterao}. In \cite{stanley2007introduction}, Stanley showed
that $\A(G)$ is free if and only if $G$ is a chordal graph. See also \cite{suyama2015vertex-weighted-a} and \cite{suyama2017signed} for more general results.

We verify this result for a given graphical arrangement in CoCoA.

\begin{lstlisting}[alsolanguage=cocoa]
/**/ use S::=QQ[x,y,z,w];
/**/ G:=[[1,2],[1,3],[1,4],[2,4],[3,4]];
/**/ A:=ArrGraphical(S,G);
/**/ ArrDerMod(A);
matrix( /*RingWithID(18935, "QQ[x,y,z,w]")*/
 [[1, 0, 0, 0],
  [1, x-y, 0, 0],
  [1, x-z, x*z-z^2-x*w+z*w, x*y-y*z-x*w+z*w],
  [1, x-w, 0, x*y-x*w-y*w+w^2]])
/**/ ArrExponents(A);
[0,  1,  2,  2]
/**/ B:=ArrDeletion(A,3);
/**/ IsArrFree(B);
false
\end{lstlisting}

\section{Algebras}

In \cite{orlik1994commutative}, Orlik and Terao introduced a commutative analogue of the Orlik-Solomon algebra in order to answer a question of Aomoto
related to cohomology groups of a certain ``twisted'' de Rham chain complex. 
The crucial difference between the Orlik-Solomon algebra  and Orlik-Terao algebra is not the difference between the exterior algebra and symmetric algebra, 
but rather the fact that the Orlik-Terao algebra actually records the ``weights'' of the dependencies among the hyperplanes.

Let $\A =\{H_1, \dots, H_n\}$ be an arrangement in $V$ and $\Lambda\subseteq \{1,\dots,n\}$. If $\bigcap_{i\in\Lambda}H_i\ne\emptyset$ and $\codim(\bigcap_{i\in\Lambda}H_i)<|\Lambda|$,
then we say that $\Lambda$ is $\textbf{dependent}$. If $\Lambda$ is dependent, then there exist $c_i\in K$ such that $\sum_{i\in\Lambda}c_i\alpha_i=0$.

\begin{Definition} Let $R$ be the ring $K[y_1,\dots, y_n]$. For each dependent set $\Lambda=\{i_1,\dots, i_k\}$, let $r_\Lambda=\sum_{j=1}^kc_{i_j}y_{i_j} \in R$. Define now
$$f_\Lambda=\partial(r_\Lambda)=\sum_{j=1}^kc_{i_j}(y_{i_1}\cdots \hat{y}_{i_j}\cdots y_{i_k}),$$ and let $I$ be the ideal of $R$ generated by the $f_\Lambda$.
This ideal is called the \textbf{Orlik-Terao ideal} of $\A$.
The \textbf{Orlik-Terao algebra} $\OT(\A)$ is the quotient $R/I$. The \textbf{Artinian Orlik-Terao algebra} $\AOT(\A)$ is the quotient of 
$\OT(\A)$ by the square of the variables.
\end{Definition}
These algebras and their Betti diagrams give us a lot of information on the given arrangement, for example about its formality. See for example \cite{schenck2009orlik}.

We can construct the Orlik-Terao ideal, its Artinian version and the Betti diagram of the Orlik-Terao algebra of the Braid arrangement in CoCoA as follows:

\begin{lstlisting}[alsolanguage=cocoa]
/**/ use S ::= QQ[x,y,z];
/**/ A:=ArrBraid(S,3);
/**/ OT_A := OrlikTeraoIdeal(A); OT_A;
ideal(y[1]*y[2]-y[1]*y[3]+y[2]*y[3])
/**/ PrintBettiDiagram(RingOf(OT_A)/OT_A);
        0    1
---------------
  0:    1    -
  1:    -    1
---------------
Tot:    1    1
/**/ ArtinianOrlikTeraoIdeal(A);
ideal(y[1]*y[2]-y[1]*y[3]+y[2]*y[3],  y[1]^2,  y[2]^2,  y[3]^2)
\end{lstlisting}

In \cite{abe2018solomon}, the authors introduced a new algebra associated to a hyperplane arrangement.  This algebra can be considered as a generalization of 
the coinvariant algebras in the setting of hyperplane arrangements and it contains the cohomology rings of regular nilpotent Hessenberg varieties.

\begin{Definition} Let $\A$ be an arrangement in $V$ and $f\in S$. Then the ideal $\mathfrak{a}(\A,f)=\{\delta(f)~|~\delta\in D(\A)\}$ is called the \textbf{Solomon-Terao ideal}
of $\A$ with respect to $f$. The \textbf{Solomon-Terao algebra} of $\A$ with respect to $f$ is the quotient $\ST(\A,f)=S/\mathfrak{a}(\A,f)$.
\end{Definition}

We can construct the Solomon-Terao ideal of the reflection arrangement of type $D$ with respect to
$f$ the sum of the square of the variables in the following way in CoCoA:
\begin{lstlisting}[alsolanguage=cocoa]
/**/ use S ::= QQ[x,y,z];
/**/ A:=ArrTypeD(S,3);
/**/ f:=x^2+y^2+z^2;
/**/ SolomonTeraoIdeal(A,f);
ideal(2*x^2+2*y^2+2*z^2, 6*x*y*z, 2*x^2*y^2-2*y^4+2*x^2*z^2-2*z^4)
\end{lstlisting}

\section{Multiarrangements of hyperplanes}

A \textbf{multiarrangement} is a pair $(\A, m)$ of an arrangement $\A$ with a map $m\colon\A\to\mathbb{Z}_{\ge0}$, called the \textbf{multiplicity}.
An arrangement $\A$ can be identified with a multiarrangement with constant multiplicity $m\equiv 1$, which is sometimes called a simple arrangement. 
Define $Q(\A,m) =  \prod_{i=1}^n\alpha_i^{m(H_i)}$ and $|m| =  \sum_{i=1}^nm(H_1)$.
With this notation, the main object is the \textbf{module of vector fields logarithmic tangent} to $\A$ with multiplicity $m$ (logarithmic vector field) 
defined by 
$$D(\A,m) = \{\delta\in \Der_{V}~|~ \delta(\alpha_i) \in \ideal{\alpha_i}^{m(H_i)} S, \forall i\}.$$

The module $D(\A,m)$ is a graded $S$-module. In general, contrarily to the case of simple arrangements, we have that
$D(\A,m)$ does not coincide with $\{\delta\in \Der_{V}~|~ \delta(Q(\A)) \in \ideal{ Q(\A,m)} S\}$. 

\begin{Definition} 
Let $\A$ be central arrangement. The multiarrangement $(\A,m)$ is said to be \textbf{free with exponents $(e_1,\dots,e_l)$} 
if and only if $D(\A,m)$ is a free $S$-module and there exists a basis $\delta_1,\dots,\delta_l \in D(\A,m)$ 
such that $\pdeg(\delta_i) = e_i$, or equivalently $D(\A,m)\cong\bigoplus_{i=1}^lS(-e_i)$.
\end{Definition}

As for simple arrangements, if  $\delta_1,\dots,\delta_l \in D(\A,m)$, then $\det(\delta_i(x_j))$ is divisible by $Q(\A,m)$.
Moreover, we can generalize Theorem \ref{theo:saitocrit}.

\begin{Theorem}[Generalized Saito's criterion]
Let $\delta_1, \dots, \delta_l \in D(\A)$. Then the following facts are equivalent
\begin{enumerate}
\item $D(\A,m)$ is free with basis $\delta_1, \dots, \delta_l$, i. e. $D(\A,m) = S\cdot\delta_1\oplus \cdots \oplus S \cdot\delta_l$.
\item $\det(\delta_i(x_j))=c Q(\A,m)$, where $c\in K\setminus\{0\}$.
\item $\delta_1, \dots, \delta_l$ are linearly independent over $S$ and $\sum_{i=1}^l\pdeg(\delta_i)=|m|$.
\end{enumerate}
\end{Theorem}

Given a simple arrangement $\A$ and $H$ one of its hyperplane, we can naturally define the \textbf{Ziegler's multirestriction} (see \cite{ziegler1986multiarrangements})
as the multiarrangement $(\A^H, m^H)$, where the function $m^H\colon A^H\to\mathbb{Z}_{>0}$ is defined by
$$X\in\A^H\mapsto\#\{H'\in\A~|~H\supset X\}-1.$$
\begin{Theorem}[\cite{ziegler1986multiarrangements}] Let $\A$ be a central arrangement. If $\A$ is free with exponents $(1,e_2,\dots,e_l)$, 
then $(\A^{H_1}, m^{H_1})$ is free with exponents $(e_2,\dots,e_l)$.
\end{Theorem}
In general, the converse of the previous theorem is false. However, we have the following

\begin{Theorem}[\cite{yoshinaga2004characterization}] Assume $l\ge 4$. Then a central arrangement $\A$ is free with exponents $(1,e_2,\dots,e_l)$ if and only if the following conditions are satisfied.
\begin{enumerate}
\item $\A$ is locally free along $H_1$, i.e. $\A_X$ is free for any $X\in L(\A)$ with $X\subset H_1$ and $X\ne\emptyset$,
\item the Ziegler's multirestriction $(\A^{H_1}, m^{H_1})$ is a free multiarrangement with exponents $(e_2,\dots,e_l)$.
\end{enumerate}
\end{Theorem}

We can construct the Ziegler's multirestriction of a given arrangement and verify the previous statements in CoCoA as follows:
\begin{lstlisting}[alsolanguage=cocoa]
/**/ use S ::= QQ[x,y,z];
/**/ A:=[x,y,z,x-y,x-y-z,x-y+2*z];
/**/ A_1:=MultiArrRestrictionZiegler(A,z);A_1;
[[y[1],  1],  [y[2],  1],  [y[1]-y[2],  3]]
/**/ MultiArrDerMod(A_1);
matrix( /*RingWithID(18, "QQ[y[1],y[2]]")*/
 [[y[1]*y[2], y[1]^3],
  [y[1]*y[2], 3*y[1]^2*y[2]-3*y[1]*y[2]^2+y[2]^3]])
/**/ MultiArrExponents(A_1);
[2,  3]
/**/ ArrExponents(A);
[1,  2,  3]
\end{lstlisting}


\bibliography{bibliothesis}{}
\bibliographystyle{plain}

\end{document}